\documentclass[11pt]{article}
\usepackage{amsfonts}
\usepackage{color,xcolor}
\usepackage{amsmath}
\usepackage{latexsym,amsfonts,amssymb}
\usepackage{mathrsfs}

\newcommand{\vsss}{\vspace{.05in}}

\begin{document}

\title{Asymptotic behavior of a free boundary problem
 for the growth of multi-layer tumors in necrotic phase}

\author{Junde Wu }

\date{   Department of Mathematics, Soochow University,
  Suzhou,  \\ [0.1 cm]
  Jiangsu 215006, PR China. E-mail:\,wujund@suda.edu.cn}

\medskip

\maketitle

\begin{abstract}
  In this paper we study a free boundary problem for
  the growth of multi-layer tumors in necrotic phase.
  The tumor region is strip-like and divided into necrotic
  region and proliferating region with two free boundaries.
  The upper free boundary is tumor surface and governed 
  by a Stefan condition.  The lower
  free boundary is the interface separating necrotic region
  from proliferating region,  its evolution is implicit and 
  intrinsically governed by an obstacle problem.
  We prove that the problem has a unique flat stationary solution,
  and there exists a positive constant $\gamma_*$, such that 
  the flat stationary solution
  is asymptotically stable for  
  cell-to-cell adhesiveness $\gamma>\gamma_*$,
  and unstable for $0<\gamma<\gamma_*$.

\medskip

  {\bf 2010 Mathematics Subject Classification}:  35B40; 35K55; 35Q92; 35R35

\medskip

  {\bf Keywords}:  free boundary problem; asymptotic behavior;
  necrotic tumor

\medskip

\end{abstract}

\section{Introduction}
\setcounter{equation}{0}
\hskip 1em

  In this paper, we study the following free boundary problem modeling
  tumor growth in necrotic phase:
\begin{equation}\label{1.1}
\left\{
\begin{array}{rll}
   \Delta\sigma\,&=\sigma \chi_{\Omega^+(t)}\qquad \qquad  & \mbox{in}\;\;
   \Omega^+(t)\cup\Omega^-(t),\; t>0,
   \\ [0.2 cm]
  \Delta p\,&=-\mu(\sigma-\tilde{\sigma})  \chi_{\Omega^+(t)}+ \nu  \chi_{\Omega^-(t)}
  \quad  \qquad  & \mbox {in}\;\;
   \Omega^+(t)\cup\Omega^-(t), \; t>0,
   \\ [0.2 cm]
  \sigma\,&=\bar{\sigma},  \quad p=\gamma\kappa
  \qquad & \mbox{on}\;\; \Gamma^+(t), \; t>0,
\\ [0.2 cm]
   \sigma\,&=\hat\sigma,\quad [\![\partial_n\sigma]\!]=0 & \hbox{on}\;\; \Gamma^-(t), \; t>0,
\\ [0.2 cm]
  [\![ p ]\!]\,&=0,\quad \, [\![\partial_n p]\!]=0\quad
   & \hbox{on}\;\; \Gamma^-(t), \; t>0,
\\ [0.2 cm]
  \partial_y \sigma\,&=0, \quad \,\partial_y p=0 & \hbox{on}\;\; \Gamma_0, \; t>0,
\\ [0.2 cm]
     V\,&=-\partial_{n} p \qquad \qquad \qquad & \hbox{on}\;\; \Gamma^+(t), \; t>0,
\\ [0.2 cm]
  \Gamma^+(0)\,&=\Gamma^+_0  & \hbox{at}\;\;  t=0,
\end{array}
\right.
\end{equation}
  where $\sigma=\sigma(x,y,t)$ and $p=p(x,y,t)$ are
  unknown functions representing concentration of nutrient
  and internal pressure
  within tumor, respectively,  $\Omega^+(t)$ and $\Omega^-(t)$ are 
  unknown domains occupied by tumor proliferating cells and necrotic cells
  at time $t>0$,
  respectively,  and
$$
 \Omega^+(t):=\{(x,y)\in \Bbb R^{n-1}\times \Bbb R:\eta(x,t)<y<\rho(x,t)\},
$$
$$
 \Omega^-(t):=\{(x,y)\in \Bbb R^{n-1}\times \Bbb R: 0<y<\eta(x,t)\},
$$
  where  $\eta(x,t)$ and $\rho(x,t)$ are unknown functions satisfying 
  $0<\eta(x,t)<\rho(x,t)$ for  $x\in \Bbb R^{n-1}$ and $t>0$,
  $\Gamma^+(t)$ and $\Gamma^-(t)$ are free boundaries, 
  and
$$
  \Gamma^+(t):={\rm graph}(\rho(x, t))=
  \{(x,y)\in\Bbb R^{n-1}\times \Bbb R: y=\rho(x,t)\},
$$
$$
  \Gamma^-(t):={\rm graph}(\eta(x,t))=\{
  (x,y)\in\Bbb R^{n-1}\times \Bbb R: y=\eta(x,t)\},
$$  
  $\Gamma_0:={\rm graph}(0)$ is the fixed bottom boundary,
  $\kappa$ is the mean curvature and $V$ 
  is the outward normal velocity of 
  upper tumor surface $\Gamma^+(t)$, respectively, $\partial_n$ denotes
  the outward normal derivative  with respect to $\Omega^{+}(t)$,
  $\bar{\sigma}$, $\tilde{\sigma}$, $\hat\sigma$, $\mu$, $\nu$ and 
  $\gamma$ are positive constants, where $\bar{\sigma}$ represents constant
  external nutrient supply, $\tilde{\sigma}$ is a critical value for
  the balance of cell apoptosis and mitosis, $\hat\sigma$ is the nutrient 
  level for tumor cell necrosis, $\mu$ is the proliferation rate of tumor 
  proliferating cells,
  $\nu$ is the  dissolution rate of necrotic cells,
  and $\gamma$ is cell-to-cell adhesiveness. We assume
  that $0<\hat\sigma<\tilde\sigma<\bar\sigma$. $\chi_{\Omega^\pm(t)}$ is the
  characteristic function of $\Omega^\pm(t)$, respectively. 
  The notation $[\![p]\!]$ 
  denotes the jump of $p$ across $\Gamma^-(t)$, and
$$
 [\![p]\!]:=\Upsilon p^+-\Upsilon p^-  \qquad\mbox{for}\quad
  p^+=p\big|_{\Omega^+(t)}\;\hbox{\; and\; }\;
  p^-=p\big|_{\Omega^-(t)},
$$
  where $\Upsilon$ is the trace operator on $\Gamma^-(t)$.
  Similarly,  $[\![\partial_n p]\!]$ and $[\![\partial_n \sigma]\!]$
  denote the jump of the normal derivatives of $p$ and $\sigma$ 
  across $\Gamma^-(t)$, respectively.

  Problem (\ref{1.1}) is originated from mathematical model proposed by
  Byrne and Chaplain [\ref{byr-cha-96}], for the growth of necrotic 
  tumor in vitro which is cultivated on an impermeable support membrane, 
  and  tumor cells are multilayered. 
  The first equation describes the diffusion and consumption of nutrient 
  in tumor region; the second equation is based on Darcy's law and mass 
  conservation law;
  the third equation means constant nutrient supply and pressure is continuous
  across the upper tumor surface, by taking cell-to-cell adhesiveness into account;
  the later three lines of equations mean that nutrient concentration, pressure and 
  their normal derivatives are continuous across the lower free boundary,
  nutrient and tumor cells can not pass 
  through the bottom boundary. For 
  more details we refer to [\ref{byr-cha-96}].
  
  In the non-necrotic case, i.e., $\Omega_-(t)=\emptyset$, the corresponding
  problem of (\ref{1.1}) has been well studied. Cui and Escher [\ref{cui-esc-09}] 
  established local well-posedness and asymptotic stability of the unique flat 
  equilibrium (independent of $x$). 
  Zhou et al. [\ref{zhou-esc-cui}]  proved that there exist infinitely many bifurcation 
  stationary solutions.  It is worthy to mention that
  another extensively studied model is solid tumor spheroid model,  where 
  tumor region is sphere-shaped. For similar solid tumor spheroid models, 
  many illuminative results such as global
  well-posedness, existence of bifurcation stationary solutions and Hopf
  bifurcations, and asymptotical stability of radially symmetric equilibrium
  have been established, we refer to [\ref{cui-16},
  \ref{cui-esc-07}, \ref{cui-esc-08}, \ref{fri-07}, \ref{fri-hu-06-0}, \ref{fri-hu-06},
  \ref{fri-hu-08}, \ref{fri-rei-99}, \ref{hao-12}, \ref{wu-16},
  \ref{wu-cui-10}, \ref{wu-zhou-17}] 
  and references therein.

  In the necrotic case,  we observe that problem (\ref{1.1}) has two free boundaries.
  The evolution of the upper free boundary $\Gamma^+(t)$ is  
  governed by the equation $V=-\partial_n p$, but the evolution of the 
  lower free boundary $\Gamma^-(t)$ is implicit. This is a remarkable
  feature and make the analysis of problem (1.1) in high dimension 
  is much more difficult than the non-necrotic case. 
  By the maximum principle,  since $0<\hat\sigma<\bar\sigma$, 
  we have $\sigma(x,y,t)\equiv \hat\sigma$ in $\Omega^-(t)$, 
  and $\sigma(x,y,t)> \hat\sigma$ in $\Omega^+(t)$ at each 
  $t>0$. Let $\Omega(t):= \{(x,y)\in\Bbb R^{n-1}\times 
  \Bbb R: 0<y<\rho(x,t)\}$ for a given function $\rho(x, t)$, 
  we can rewrite $\Omega^+(t)=\{(x,y)
  \in\Omega(t): \sigma(x,y,t) >\hat\sigma\}$, 
  $\Omega^-(t)=\mbox{int}\{ (x,y)\in \Omega(t):
  \sigma(x,y,t)=\hat\sigma\}$ and 
  $\Gamma^-(t)=\partial\Omega^+(t)\cap
  \partial \Omega^-(t)$.  We see that  
  $(\sigma(x,y,t), \Gamma^-(t))$ satisfies 
   an obstacle problem:
\begin{equation}\label{1.2}
  \left\{
  \begin{array}{l}
  -\Delta \sigma +\sigma\ge0,\qquad \sigma \ge \hat \sigma, \qquad
  (-\Delta\sigma+\sigma)(\sigma-\hat\sigma)=0\qquad\mbox{in }\;\; \Omega(t),
  \\ [0.2 cm]
  \sigma=\bar\sigma \qquad
  \mbox{on }\;\; \Gamma^+(t),
  \qquad\qquad
  \partial_y\sigma=0 \qquad \mbox{on }\;\; \Gamma_0.
 \end{array}
  \right.
\end{equation}
  The regularity of free boundary of obstacle problems in high dimension is 
  very difficult to study (cf. [\ref{caf-77}, \ref{kin-sta}]).  Even for
  smooth domain $\Omega(t)$, the solution $\sigma\notin C^2(\Omega(t))$.
  It makes a main difficulty arising 
  in necrotic tumor model from non-necrotic case.   To solve this
  problem, 
  we first show there is a unique flat stationary solution 
  (see Section 2). It implies that we can consider obstacle problem (\ref{1.2}) 
  in a small neighborhood of the flat stationary solution.  Motivated by
  Cui [\ref{cui-16}] and Hamilton [\ref{ham-82}],
  by using Nash-Moser implicit function theorem,
  for given $\Gamma^+(t)={\rm graph}(\rho(x,t))$ closed to 
  the flat equilibrium, we prove that
   the solution $(\sigma(x,y,t), \Gamma^-(t))$ 
  is smoothly depending on $\rho(x,t)$, and $\Gamma^-(t)$ 
  is actually smooth in space variables (see Lemma 3.2). 
  Then we further solve the first six lines of equations of problem $(\ref{1.1})$ and 
  get the solution $p(x,y,t)$ which also smoothly depends on $\rho(x,t)$, 
  finally by the second last equation
  $V=-\partial_n p$ we reduce problem (\ref{1.1}) into an abstract differential
  equation $\partial_t\rho+\Psi(\rho)=0$ only containing
   function $\rho(x,t)$.  In suitable Banach spaces, we show this abstract 
  differential equation is of parabolic type and the local well-posedness follows
  by geometric theory of parabolic differential equations. By a delicate analysis and 
  computation, we study the spectrum of the linearized operator at the
  flat stationary solution,  and by linearized stability principle we can 
  get asymptotic stability of the flat stationary solution. 
  
  To give a precise statement
  of our main results, we introduce some notations.
  
  In this paper, we only consider the case $n=2$, and
  the higher-dimensional case can be treated similarly. 
  We denote the solution of problem (\ref{1.1}) by
  $(\sigma, p, \eta, \rho)$, with 
  $\Gamma^-(t)={\rm graph}(\eta(x,t))$
  and $\Gamma^+(t)={\rm graph}(\rho(x,t))$. 
  For the sake of simplicity, we impose that 
\begin{equation}\label{1.3}
  \sigma(x,y,t),\; p(x,y,t), \; \eta(x,t),\; \rho(x,t) \mbox{ are } 
  2\pi\mbox{-periodic  in } x\in\mathbb{R}.
\end{equation}
  We identify $\mathbb {S}=\mathbb{R}/{2\pi}\mathbb{Z}$, and
  identify continuous
  2$\pi$-periodic function space $C_{per}(\mathbb{R})=C(\mathbb{S})$. 
  Given $s>0$, 
  we denote by $BUC^{s}(\Bbb S)$ the space of all bounded and uniformly 
  H\"older continuous functions on $\Bbb S$ of order $s>0$. 
  Let $h^{s}(\Bbb S)$ denote the
  little H\"older space, a closure of $BUC^\infty(\Bbb S)$
  in $BUC^{s}(\Bbb S)$.  Similarly, we denote by $h^{s}(\Omega)$
  the closure of $BUC^{\infty}(\Omega)$ in $BUC^{s}(\Omega)$ for
  bounded open domain $\Omega$ in $\Bbb R^2$.

  Our first main result is stated as follows:
\medskip

  {\bf Theorem 1.1} \ \ {\em  Let $\tilde\sigma>\hat\sigma>0$ be given.
  There exists a positive constant $\sigma_*>\tilde\sigma$
  depending only on $\hat\sigma$ and $\tilde\sigma$, such that 
  free boundary problem $(1.1)$ has a unique flat stationary solution 
  $(\sigma_s, p_s, \eta_s, \rho_s)$ if and only if $\bar{\sigma}>
  \sigma_*$.
}

\medskip  
  
  We shall prove this result in Section 2.  Recall that in non-necrotic case
   (see Theorem 2 of [\ref{cui-esc-09}]), there exists a unique 
   non-necrotic flat stationary solution for all $\bar\sigma>\tilde\sigma$.
   It is an interesting difference that necrotic flat stationary solution
   does not exist for $\tilde\sigma<\bar\sigma\le \sigma_*$.

\medskip   

  Our second main result is about asymptotic stability of the flat stationary solution.
     
\medskip

  {\bf Theorem 1.2} \ \ {\em 
  $(i)$ There exists a positive threshold value $\gamma_*$
  of cell-to-cell adhesiveness such that for any $\gamma>\gamma_*$,
  the flat stationary solution $(\sigma_s, p_s, \eta_s, \rho_s)$ is asymptotically 
  stable in the following sense: There exists a constant $\epsilon>0$ such that
  if $\rho_0\in h^{4+\alpha}(\mathbb{S})$,
  $\|\rho_0\|_{h^{4+\alpha}(\mathbb{S})}<\epsilon$ and $\Gamma^+_0=
  {\rm graph}(\rho_s+\rho_0)$,
  then
  the solution $(\sigma, p, \eta, \rho)$ of problem $(\ref{1.1})$ exists for 
  all $t > 0$ and converges to $(\sigma_s, p_s, \eta_s, \rho_s)$ 
  exponentially fast as  $t\to+\infty$.

  $(ii)$
  If  $0<\gamma<\gamma_*$, then $(\sigma_s, p_s, \eta_s, \rho_s)$
  is unstable.}
\medskip
  
  The above result implies that cell-to-cell
  adhesiveness $\gamma$ plays an important role on tumor's
  stability. A smaller value of $\gamma$ may make tumor
  more aggressive. The threshold value $\gamma_*$ of
  cell-to-cell adhesiveness is given by (\ref{4.19})
  and (\ref{4.21}), and $\gamma_*$ can be regarded as a
  function of the dissolution rate $\nu$. By
  $\displaystyle {d \gamma_*/ d \nu}\le0$, 
  we see that a 
  smaller value of $\nu$ may make
  tumor more unaggressive.  While in the limiting case $\nu=0$,
   the flat stationary solution is not asymptotically stable 
   for  all $\gamma>0$.  (see Remark 5.2).

\medskip

  The structure of the rest of this paper is arranged as follows.
  In the next section, we study the existence and uniqueness 
  of flat stationary solution. In Section 3, by using implicit 
  function theorem and classical theory of elliptic equations
  we reduce free boundary
  problem (\ref{1.1}) into a Cauchy problem in little H\"older
  spaces, and establish the local well-posedness. 
  Section 4 is devoted to study the 
  linearized problem at flat stationary solution and 
  compute eigenvalues. In the
  last section we make  stability analysis and give a proof of 
  Theorem 1.2.
  
\medskip
\hskip 1em

\section{Flat stationary solution}
\setcounter{equation}{0}
\hskip 1em

  In this section, we study the existence and uniqueness
  of flat stationary solution of free boundary problem 
  (\ref{1.1}).  
  
   We denote flat stationary solution by $(\sigma_s(y), p_s(y), \eta_s, \rho_s)$ with
   $0<\eta_s<\rho_s$. It satisfies the following problem
  
\begin{equation}\label{2.1}
\left\{
\begin{array}{rll}
   \sigma_s''(y)\,&=\sigma_s(y),  \qquad \quad  
    p_s''(y)=-\mu(\sigma_s(y)-\tilde{\sigma}) \qquad \qquad & \hbox{ for }
   \eta_s<y<\rho_s,
   \\ [0.2 cm]
   \sigma_s''(y)\, &= 0, \qquad\quad\quad\;\;\; p_s''(y)=\nu
  \quad  \qquad  &\hbox{ for }\;  0<y<\eta_s,
\\ [0.2 cm]
  \sigma_s(\rho_s)\,&=\bar{\sigma},  \qquad\qquad\;\; p_s(\rho_s)=0,
\\ [0.2 cm]
   \sigma_s(\eta_s)\,&=\hat\sigma,\qquad\qquad \;\;\sigma_s'(\eta_s)=0 ,
   \\ [0.2 cm]
   p_s(\eta_s^+)\,&=p_s(\eta_s^-),\quad\quad\; p_s'(\eta_s^+)=p_s'(\eta_s^-),
\\ [0.2 cm]
   \sigma_s'(0)\,&=0, \qquad p_s'(0)=0 ,\qquad p_s'(\rho_s)=0.
\end{array}
\right.
\end{equation}
  We easily get that
\begin{equation}\label{2.2}
  \sigma_s(y)=
  \left\{
  \begin{array}{ll}
  \displaystyle {\bar\sigma \sinh (y-\eta_s)
  +\hat\sigma\sinh(\rho_s-y)\over\sinh(\rho_s-\eta_s)}
  \qquad\qquad\qquad\qquad\quad\;\; & \mbox{for}\;\;\eta_s\le y\le \rho_s,
  \\ [0.2 cm]
  \hat\sigma\qquad
  & \mbox{for}\;\;  0<y<\eta_s,
  \end{array}
  \right.
\end{equation}  
\begin{equation}\label{2.3}
  p_s(y)=
  \left\{
  \begin{array}{ll}
  \displaystyle {\mu\over 2}\tilde\sigma(y^2-\rho_s^2)
  +(\nu-\mu\tilde\sigma)(y-\rho_s)\eta_s+\mu(\bar\sigma-\sigma_s(y))
  \qquad & \mbox{for}\;\;\eta_s\le y\le \rho_s,
  \\ [0.3 cm]
  \displaystyle {\nu\over 2} (y^2-\eta_s^2)+p_0 \qquad
  & \mbox{for}\;\;  0<y<\eta_s,
  \end{array}
  \right.
\end{equation}  
where $p_0=
  {\mu\over 2}\tilde\sigma(\eta_s^2-\rho_s^2)+(\nu-\mu\tilde\sigma)
  (\eta_s-\rho_s)\eta_s
  +\mu(\bar\sigma-\hat\sigma)$.

  By $\sigma_s'(\eta_s)=0$, there holds
\begin{equation}\label{2.4}
 \cosh(\rho_s-\eta_s)={\bar\sigma\over\hat\sigma}.
\end{equation}
  Using this formula,
\begin{equation}\label{2.5}
  \sigma_s'(\rho_s)={\bar\sigma\cosh(\rho_s-\eta_s)-\hat\sigma\over
  \sinh(\rho_s-\eta_s)}
  = \sqrt{\bar\sigma^2-\hat\sigma^2}.
\end{equation}
  By $p_s'(\rho_s)=0$ we have $ \mu\tilde\sigma\rho_s+
  (\nu-\mu\tilde\sigma)\eta_s-\mu\sigma_s'(\rho_s)=0$.
  It implies that
\begin{equation}\label{2.6}
  (\nu-\mu\tilde\sigma)\eta_s+\mu\tilde\sigma\rho_s=\mu\sqrt{\bar\sigma^2-\hat\sigma^2}.
\end{equation}
  Then from (\ref{2.4}) and (\ref{2.6}), we obtain
\begin{equation}\label{2.7}
  \eta_s={\mu\over\nu}\big( \sqrt{\bar\sigma^2-\hat\sigma^2} - \tilde\sigma
  \ln (\bar\sigma+\sqrt{\bar\sigma^2-
  \hat\sigma^2}) +\tilde\sigma\ln \hat\sigma
  \big),
\end{equation}
\begin{equation}\label{2.8}
  \rho_s=\eta_s+\ln (\bar\sigma+\sqrt{\bar\sigma^2-
  \hat\sigma^2}) -\ln \hat\sigma.
\end{equation}
  Clearly,  $\rho_s>\eta_s$ if and only if $0<\hat\sigma<\bar\sigma$.
  
  Next, we only need to make sure $\eta_s>0$ for
  $0<\hat\sigma<\min\{\tilde\sigma,\bar\sigma\}$.
  Define a function 
$$
  f(a,r):=\sqrt{r^2-1}-a \ln (r+\sqrt{r^2-1})\qquad\mbox{for}\;\;
  r>1,\; a>1.
$$
  Note that
$$
  \partial_r f(a,r) = {r-a\over \sqrt{r^2-1}},\;\qquad\; f(a,a)<0\qquad\;
  \mbox{for}\;\; r>1,\; a>1.
$$
  It follows that for any $a>1$, there exists a positive constant $a_*>a$
  such that
$$
  f(a,r)\left\{
  \begin{array}{l}
  <0,\qquad 1< r<a_*,
  \\ [0.2 cm]
  =0, \qquad r=a_*,
  \\ [0.2 cm]
  >0,\qquad r>a_*.
  \end{array}
  \right.
$$
  By (\ref{2.7}) we see that $\displaystyle \eta_s = {\mu\hat\sigma\over \nu}
   f( {\tilde\sigma\over \hat\sigma}, {\bar\sigma\over\hat\sigma})$.
  Recall that $0<\hat\sigma<\min\{\tilde\sigma,\bar\sigma\}$. We immediately 
  obtain that there exists a positive constant $\sigma_*>\tilde\sigma$
  depending only on $\tilde\sigma$ and $\hat\sigma$,
  such that $\eta_s>0$ for $\bar\sigma>\sigma_*$ and $\eta\le 0$ for 
  $\bar\sigma\le \sigma_*$.
  
  In conclusion, we have
  
\medskip

{\bf Theorem 2.1}\ \ {\em Assume $0<\hat\sigma<\tilde\sigma$.
  There exists a positive constant $\sigma_*>\tilde\sigma$
  depending only on $\tilde\sigma$ and $\hat\sigma$, such that for 
  $\bar\sigma>\sigma_*$, problem $(\ref{1.1})$ has a unique flat stationary solution
  $(\sigma_s, p_s,\eta_s,\rho_s)$ given by $(\ref{2.2})$, $(\ref{2.3})$, 
  $(\ref{2.7})$ and $(\ref{2.8})$.
  If  $\bar\sigma\le\sigma_*$, problem $(\ref{1.1})$ has no flat stationary solution.
  }

\medskip

   It is interesting to compare this result with non-necrotic case.
   From Theorem 2 of [\ref{cui-esc-09}], we see that there exists a unique 
   non-necrotic flat stationary solution for $0<\tilde\sigma<\bar\sigma$.
   But in necrotic case, we see for $\tilde\sigma<\bar\sigma\le \sigma_*$,
   necrotic flat stationary solution does not exist.

\medskip
\hskip 1em

\section{Reduction and Well-posedness}
\setcounter{equation}{0}
\hskip 1em

  In this section, we  reduce free boundary problem (\ref{1.1})
  into a Cauchy problem in little H\"older spaces, and study the 
  local well-posedness.

  First, we transform free boundary problem (\ref{1.1}) 
  into an equivalent problem on a fixed domain. Later on,
  we always assume $0<\hat\sigma<\tilde\sigma<\sigma_*<\bar\sigma$.
  By Theorem 2.1, problem (\ref{1.1}) has a unique flat stationary solution 
  $(\sigma_s, p_s,\eta_s,\rho_s)$. Denote
$$
  \Omega_s=\{(x,y)\in\Bbb S\times\Bbb R: 0<y<\rho_s\},\qquad
  \Bbb D_s=\{(x,y)\in\Bbb S\times\Bbb R: 0<y<\eta_s\},
$$
$$
  \Gamma_s=\Bbb S\times\{\rho_s\},\qquad
  J_s=\Bbb S\times\{\eta_s\},\qquad
  \Gamma_0=\Bbb S\times\{0\},\qquad
  \Bbb E_s=\Omega_s\backslash \overline{\Bbb D}_s.
$$
    Let  $r_0:=(\rho_s-\eta_s)/8$,
    $\delta\in (0, r_0)$ and $\alpha\in(0,1)$, set
\begin{equation}\label{3.1}
  \mathcal{O}_\delta:=\{\rho\in h^{4+\alpha}(\Bbb S):
  \,\|\rho\|_{h^{4+\alpha}(\Bbb S)}<\delta \}.
\end{equation}
  For  $\rho,\eta\in \mathcal{O}_\delta$, we denote
$$
\begin{array}{rl}
  \Omega_\rho=\{(x,y)\in \Bbb S\times \Bbb R: 0<y<\rho_s+\rho(x)\},
\qquad
 &  \Bbb D_\eta=\{(x,y)\in \Bbb S\times \Bbb R: 0<y<\eta_s+\eta(x)\},
\\[0.2 cm]
  \Gamma_\rho=\{(x,y)\in\Bbb S\times\Bbb R: y=\rho_s+\rho(x)\},
  \qquad
  & J_\eta=\{(x,y)\in\Bbb S\times\Bbb R: y=\eta_s+\eta(x)\},
\end{array}
$$ 
$$
  \Omega_{\rho,\eta}=\Omega_\rho\backslash \overline {\Bbb D}_\eta
    \qquad\quad\mbox{and}\qquad\quad
  \Bbb E_\eta=\Omega_s\backslash \overline {\Bbb D}_\eta.
  \qquad\qquad\qquad
$$
   Choose a function $\varphi \in C^\infty(\Bbb R)$ such that
\begin{equation}\label{3.2}
  0\le\varphi(y)\le1,
\qquad
  \varphi(y)=
  \left\{
  \begin{array}{l}
  1,\quad \mbox{for}\;|y| \le \delta,
  \\ [0.2 cm]
  0,\quad \mbox{for}\;|y| \ge 3\delta,
  \end{array}
  \right.
  \qquad
  \sup|\varphi'(y)|< {1/ \delta}.
\end{equation}
   Given $\rho\in\mathcal{O}_\delta$, we introduce
   a mapping  
$$
  \Phi_\rho: \Omega_s \to \Omega_\rho,
  \qquad 
   (x,y)\to \big(x, y+\varphi(y-\rho_s)\rho(x)\big).
$$
  Clearly,  $\Phi_\rho(\Omega_s)=\Omega_\rho$,
  $\Phi_\rho(\Gamma_s)=\Gamma_\rho$
  and $\Phi_\rho$ is a $h^{4+\alpha}$ diffeomorphism
  from $\Omega_s$ onto $\Omega_\rho$.
  Moreover, for any $\eta\in\mathcal O_\delta$, 
  $\Phi_\rho$  is the identity mapping on $\Bbb D_\eta$.
  Define the induced push-forward operator $\Phi^\rho_*$, 
  and pull-back operator $\Phi^*_\rho$ by
\begin{equation}\label{3.3}
  \Phi^\rho_* u=u\circ \Phi_\rho^{-1} \qquad \mbox{for}\;\;
  u\in C(\Omega_s),
\qquad\;\;
 \Phi^*_\rho v=v \circ \Phi_\rho \qquad \mbox{for}\;\;
  v\in C(\Omega_\rho).
\end{equation}
  Next,  we introduce the following transformed operators:
\begin{equation}\label{3.4}
 \;\;\; \mathcal A(\rho)u:=\Phi_{\rho}^*\Delta(\Phi_*^{\rho}u),
  \qquad\;
  \mathcal B(\rho)u:=\langle\nabla(\Phi_*^\rho u)|_{\Gamma_\rho},
  {\bf n}_\rho\rangle
  \qquad\;\mbox{for}\;\; u\in H^2(\Omega_s),
\end{equation}
  where ${\bf n}_\rho=(-\rho_x,1)$ is the outward normal on
  $\Gamma_\rho$,  $\langle\cdot,\cdot \rangle$ denotes the
  Euclidean inner product, and
  $\displaystyle H^2(\Omega_s)$ stands for Sobolev space.  
  By Lemma 2.2 of [\ref{esc-sim-97-1}],
  we have
\begin{equation}\label{3.5}
\left\{
\begin{array}{l}
  \mathcal A\in C^\infty\big(\mathcal O_\delta,L(h^{k+2+\alpha}(\Omega_s),
 h^{k+\alpha}(\Omega_s))\big),\;\;\;\; 0\le  k\le 2,
\\ [0.2 cm]
  \mathcal B \in C^\infty\big(\mathcal O_\delta,L(h^{k+1+\alpha}(\Omega_s),
 h^{k+\alpha}(\Bbb S))\big),\quad \;\;\; 0\le k\le 3.
 \end{array}
 \right.
\end{equation}
   Denote by $\mathcal K(\rho)$ the transformed mean curvature 
   on $\Gamma_{\rho}$ and 
\begin{equation}\label{3.6}
  \displaystyle\mathcal K(\rho)=-(1+\rho_x^2)^{-{3\over2}}\rho_{xx}.
\end{equation}

   For some $T>0$, and a function 
   $\rho\in C([0,T), \mathcal O_\delta)\cap C^1([0,T),
  h^{1+\alpha}(\mathbb{S}))$, we identify  $\rho(x,t)=\rho(t)(x)$
  for $t\in [0,T)$ and $x\in\mathbb{S}$. 
  By an elementary analysis,  the outward normal velocity  $V$
  of tumor surface $\Gamma_{\rho(t)}$ can be given by 
$$
  V={\rho_t/ \sqrt{1+\rho_x^2}}.
$$

  Let $\chi_{\Bbb D_\eta}$ and $\chi_{\Bbb E_\eta}$ be the 
  characteristic functions of $\Bbb D_\eta$ and $\Bbb E_\eta$, 
  respectively.
   Rewrite 
$$
  \Gamma^+_0={\rm graph}(\rho_s+\rho_0)\qquad \mbox{for\; some }\;\;
  \rho_0\in\mathcal O_\delta,  
$$
  and
$$
  u(x,y,t)=\Phi^*_\rho\sigma (x,y, t), \quad
  v(x,y,t)=\Phi^*_\rho\, p(x,y,t).
$$
  One can easily check that  free boundary problem (\ref{1.1}) is
  transformed into the following problem:
\begin{equation}\label{3.7}
\left\{
\begin{array}{rll}
   \mathcal A(\rho)u\,&=u \chi_{\Bbb E_\eta} \qquad & \hbox{in}
   \;\;\Omega_s,\; t>0,
   \\ [0.2 cm]
  \mathcal A(\rho)v\,&=-\mu(u-\tilde{\sigma}) \chi_{\Bbb E_\eta} 
  +\nu \chi_{\Bbb D_\eta}
  \quad\;\;  \qquad  &\hbox{in}\;\;
   \Omega_s, \; t>0,
\\ [0.2 cm]
  u\,&=\bar{\sigma},  \quad v=\gamma\mathcal K(\rho)
  \qquad & \hbox{on}\;\; \Gamma_s, \; t>0,
\\ [0.2 cm]
   u\,&=\hat\sigma,\quad [\![\partial_n u]\!]=0 & \hbox{on}\;\; J_{\eta}, \; t>0,
\\ [0.2 cm]
   [\![v]\!]\,&=0,\quad \,[\![\partial_n v]\!]=0\quad
   & \hbox{on}\;\; J_{\eta}, \; t>0,
\\ [0.2 cm]
  \partial_y u\,&=0, \quad \,\partial_y v=0 & \hbox{on}\;\; \Gamma_0, \; t>0,
\\ [0.2 cm]
     \partial_t \rho\,&=-\mathcal B(\rho) v \qquad \qquad \qquad & \hbox{on}\;\; \Bbb S, 
     \;\;\, t>0,
\\ [0.2 cm]
  \rho(0)\,&=\rho_0  & \hbox{on}\;\;  \Bbb S,\;\;\, t=0.
\end{array}
\right.
\end{equation}

  By above transformation, we have

 \medskip

{\bf Lemma 3.1} \ \ {\em  A quadruple $(u, v, \eta, \rho)$
  is a solution of problem $(\ref{3.7})$ if and only if 
  the quadruple $(\sigma, p,\eta_s+\eta,\rho_s+\rho)$ 
  is a solution of problem $(\ref{1.1})$ in the neighborhood 
  of $(\sigma_s, p_s,\eta_s,\rho_s)$,   
  with $\sigma=\Phi_*^{\rho}u$ and 
  $p=\Phi_*^{\rho} v$.
}

\medskip

  Next we further reduce problem (\ref{3.7}) into a Cauchy problem 
  in little H\"older space for $\rho$ only.  
  Given  $\rho\in \mathcal O_\delta$, we consider the following 
  problem:  
\begin{equation}\label{3.8}
\left\{
\begin{array}{l}  
     \mathcal A(\rho)u=u \chi_{\Bbb E_\eta} \qquad  \hbox{in}\;\;
   \Omega_s,
   \\ [0.2 cm]
  u\big|_{\Gamma_s}=\bar{\sigma},
  \qquad  \partial_y u\big|_{\Gamma_0}=0,
\\ [0.2 cm]
   u=\hat\sigma,\quad [\![\partial_n u]\!]=0 \quad \hbox{on}\;\; J_{\eta}.
\end{array}
\right.
\end{equation}
  For any $\eta\in\mathcal O_\delta$, by the maximum principle, 
  $u\equiv \hat\sigma$ in $\Bbb D_\eta$, 
  and $u>\hat\sigma$ in $\Bbb E_\eta$. On the other hand, since
  $u=\hat\sigma$ on $J_\eta$, we have
$$
  \partial_n u={u_x\eta_x-u_y\over\sqrt {1+\eta_x^2}}
  =-u_y\sqrt{1+\eta_x^2} \qquad\quad \mbox{on}\;\;J_\eta.
$$
  It implies that $[\![\partial_n u]\!]=0$ is equivalent to 
  $\partial_y u = 0$ on $J_\eta$.
  Hence for problem (\ref{3.8}), we only need to solve 
\begin{equation}\label{3.9}
\left\{
\begin{array}{rlr}  
     \mathcal A(\rho)u\,&=u \qquad\;\;  & \hbox{in}\;\;
   \Bbb E_\eta,
   \\ [0.2 cm]
  u \, &=\bar{\sigma}\qquad & \hbox{on}\;\; \Gamma_s,
\\ [0.2 cm]
   \partial_y u\,& =0\qquad & \hbox{on}\;\; J_{\eta},
 \\ [0.2 cm]
   u\, &=\hat\sigma   \qquad & \hbox{on}\;\; J_{\eta}.
\end{array}
\right.
\end{equation}
 
   Recently, Cui [\ref{cui-16}] studied a similar obstacle 
   problem based on Nash-Moser implicit
   function theorem. 
   Motivated by this method and with some
   modifications to the proof of Theorem 5.2 of
   [\ref{cui-16}] ,  
   we have the following result:
   
\medskip

  {\bf Lemma 3.2} \ \ {\em There exists a constant $\delta_1\in(0,r_0)$, 
  such that
  for any $\rho\in\mathcal O_{\delta_1}$, problem $(\ref{3.9})$ has a unique solution
  $(u,\eta)$ satisfying  $u\in h^{4+\alpha}(\Bbb E_\eta)$ and 
  $\eta\in C^\infty(\Bbb S)$. Moreover, the mapping $\rho\mapsto (u,\eta)$
  from $\mathcal O_{\delta_1}$ to $h^{4+\alpha}(\Bbb E_\eta)\times C^\infty(\Bbb S)$
  is smooth.}
  
\medskip

  {\bf Proof.} \ \      Denote 
$$ 
   \Bbb H:=\{(x,y)\in\Bbb S\times\Bbb R:
   {\rho_s+\eta_s\over 2}<y<\rho_s\},
   \qquad
   \Bbb K:= \{(x,y)\in\Bbb S\times\Bbb R:
   \eta_s<y<{\rho_s+\eta_s\over 2}\}.
$$
  Let $r_1:=\min\{(\rho_s-\eta_s)/8,\eta_s/8\}\le r_0$ and $\delta\in(0,r_1)$. 
  For any $m\in\Bbb N$, $m\ge4$ and $\alpha\in(0,1)$, denote
$$
  \widetilde{\mathcal O}_\delta^{m+\alpha}:=\{\eta\in h^{m+\alpha}(\Bbb S):
  \|\eta\|_{h^{4+\alpha}(\Bbb S)}<\delta\}.
$$
  For
  any $\eta\in\widetilde{\mathcal O}_\delta^{m+\alpha}$, we introduce
   a mapping  
$$
  \widetilde\Phi_\eta: \Bbb R^2 \to \Bbb R^2,
  \qquad 
   (x,y)\to \big(x, y+\varphi(y-\eta_s)\eta(x)\big),
$$
  where $\varphi$ is a smooth function given by (\ref{3.2}).
  We have 
  $\widetilde\Phi_\eta$ is a $h^{m+\alpha}$ diffeomorphism
  from $\Bbb E_s$ onto $\Bbb E_\eta$, and 
  $\widetilde\Phi_\eta$ is the identity mapping on $\Bbb H$.
  Similarly as (\ref{3.3}), we can define the 
  push-forward operator $\widetilde\Phi^\eta_*$, 
  and pull-back operator $\widetilde\Phi^*_\eta$
  induced by $\widetilde\Phi_\eta$, and for any 
  $\rho\in\mathcal O_\delta$ and 
  $\eta\in\widetilde{\mathcal O}_\delta^{m+\alpha}$,
  we define an operator
$$
  \mathscr A(\rho,\eta)v:=
  \widetilde\Phi_{\eta}^*\mathcal A(\rho)
  (\widetilde\Phi_*^{\eta}v)
  \qquad\mbox{for}\;\; v\in BUC^2(\Bbb E_s).
$$ 
  Notice that for any $\rho\in\mathcal O_\delta$ and 
  $\eta\in\widetilde{\mathcal O}_\delta
  ^{m+\alpha}$, 
   $\mathcal A(\rho)\equiv\Delta $ in $ \Bbb K$, 
   we see that $\mathscr A(\rho,\eta)$ is independent of $\eta$ on $\Bbb H$,
   and independent of $\rho$ on $\Bbb K$. Moreover,
   $\mathscr A(\rho,\eta)$ is uniformly elliptic and 
   by Lemma 2.2 in [\ref{esc-sim-97-1}], we have  
$$
  \mathscr A\in C^\infty\big(\mathcal O_\delta\times 
  \widetilde{\mathcal O}_\delta
  ^{m+\alpha},L(h^{4+\alpha}(\Bbb E_s)\cap h^{m+\alpha}(\Bbb K),
  h^{2+\alpha}(\Bbb E_s)\cap  h^{m-2+\alpha}(\Bbb K))\big).
$$   
  Set $\tilde u =u\circ \widetilde\Phi_\eta$. The first three 
  equations of $(\ref{3.9})$ is equivalent to 
\begin{equation}\label{3.10}
     \mathscr A(\rho,\eta)\tilde u=\tilde u \quad  \hbox{in}\;\;
   \Bbb E_s,
   \qquad\;
  \tilde u =\bar{\sigma}\quad \hbox{on}\;\; \Gamma_s,
  \qquad\;
  \partial_y \tilde u =0 \quad \hbox{on}\;\; J_{s}.
\end{equation}  
  By well-known regularity theory of second-order elliptic differential equations, 
  problem $(\ref{3.10})$ has a unique solution $\tilde u:=\widetilde{\mathcal U}(\rho,\eta)
 \in h^{4+\alpha}(\Bbb E_s)$, and 
  by Lemma 2.3 in [\ref{esc-sim-97-1}], for $m\ge 4$, 
\begin{equation}\label{3.11}
  \widetilde{\mathcal U}\in C^\infty(\mathcal O_\delta\times \widetilde{\mathcal O}_\delta^{m+\alpha}, 
  h^{4+\alpha}(\Bbb E_s)).
\end{equation}
   Next, we further show some much more profound properties of  
   $\;\widetilde{\mathcal U}$. 
   Recall from Part II.1 of  Hamilton [\ref{ham-82}],  Banach space 
   $h^{4+\alpha}(\Bbb S)$ can be regarded as a tame Fr\'echet space,  
   $C^\infty(\Bbb S)$ with a collection of seminorms $\{\| \;\; \|_{h^{m+\alpha}(\Bbb S)}, 
   m=0,1,2\cdots\}$, and $BUC^\infty(\Bbb K)$ with a collection of seminorms 
   $\{\| \;\; \|_{h^{m+\alpha}(\Bbb K)}, m=0,1,2\cdots\}$ are both
   tame Fr\'echet spaces.  Denote 
$$
  \widetilde{\mathcal O}_\delta^\infty =\{\eta\in C^\infty(\Bbb S): 
  \|\eta\|_{h^{4+\alpha}(\Bbb S)}<\delta\}.
$$
    By Theorem 3.3.5 
    in Part II of [\ref{ham-82}],  we have
\begin{equation}\label{3.12}
  \widetilde{\mathcal U} \; \mbox{  is a smooth tame mapping from }
   \; \mathcal O_\delta\times
   \widetilde{\mathcal O}_\delta^\infty \; \mbox{ to } \;
    h^{4+\alpha}(\Bbb E_s)\cap BUC^\infty (\Bbb K),
\end{equation}
     which means that for all $m\in \Bbb N$, $m\ge 4$,  
\begin{equation}\label{3.13}
  \widetilde{\mathcal U}\in C^\infty(\mathcal O_\delta\times 
  \widetilde{\mathcal O}_\delta^{m+\alpha},  
  h^{4+\alpha}(\Bbb E_s)\cap h^{m+\alpha}(\Bbb K)),
\end{equation}
  and
\begin{equation}\label{3.14}
  \| \widetilde{\mathcal U}(\rho,\eta)\|_{h^{4+\alpha}(\Bbb E_s)}+
  \| \widetilde{\mathcal U}(\rho,\eta)\|_{h^{m+\alpha}(\Bbb K)}  \le C_m (1+
  \|\rho\|_{h^{4+\alpha}(\Bbb S)}+\|\eta\|_{h^{m+\alpha}(\Bbb S)}),
\end{equation}
  where $C_m$ is a positive constant dependent on $m$.  

  We define a mapping $F:\mathcal O_\delta\times\widetilde{\mathcal O}_\delta^{m+\alpha}
  \to h^{m+\alpha}(\Bbb S)$ by
$$
  F(\rho,\eta)= \widetilde{\mathcal U}(\rho,\eta) \Big|_{J_s}-\hat\sigma.
$$
  It is easy to see that $F\in C^\infty(\mathcal O_\delta\times
  \widetilde{\mathcal O}_\delta^{m+\alpha}, h^{m+\alpha}(\Bbb S))$.
  Moreover, by (\ref{3.12}) we have
\begin{equation}\label{3.15}
  F \; \mbox{  is a smooth tame mapping from }
   \; \mathcal O_\delta\times
   \widetilde{\mathcal O}_\delta^\infty \; \mbox{ to } \;
    C^\infty (\Bbb S).
\end{equation}
   Clearly, $F(0,0)=0$ and problem (\ref{3.9}) is equivalent to the equation
   $F(\rho,\eta)=0$.

  Next we compute the Fr\'echet derivative of $F$ with respect to $\eta$ 
  at  $(\rho,\eta)\in \mathcal O_\delta\times \widetilde{\mathcal O}_\delta^\infty$, 
  which is denoted  by $D_\eta F(\rho,\eta)$. 
  Let $\mathcal U(\rho,\eta)$ be the solution of the first three equations of problem (\ref{3.9}).
  For any $\zeta\in C^\infty(\Bbb S)$, we easily verify that 
$$
  D_\eta F(\rho,\eta)\zeta= \mathcal Z(\rho,\eta, \zeta)\big|_{J_\eta},
$$
  where 
 $z=\mathcal Z(\rho,\eta,\zeta)$ 
 is the solution of the following problem
\begin{equation}\label{3.16}
   \mathcal A(\rho) z=z \quad   \mbox{in}\;\;
   \Bbb E_\eta,
   \qquad
  z=0 \quad  \mbox{on}\;\; \Gamma_s,
 \qquad
 \partial_ y z=-\partial_{yy}\, \mathcal U(\rho,\eta)  \zeta
  \quad  \mbox{on}\;\; J_\eta.
\end{equation}
  Since $\mathcal U(0,0)=\sigma_s\big|_{\Bbb E_s}$, we have 
  $\partial_{yy}\,\mathcal U(0,0)\big|_{J_s}=\sigma_s''(\eta_s^+)=\hat\sigma>0$.
  Thus for sufficiently small $\delta>0$, we have 
  $\partial_{yy}\,\mathcal U(\rho,\eta)\big|_{J_\eta}>\hat\sigma/2$
  for $(\rho,\eta)\in \mathcal O_\delta\times \widetilde{\mathcal O}_\delta^\infty$.  
  By (\ref{3.16}),  for any $\xi\in C^\infty(\Bbb S)$ we have
$$
  [D_\eta F(\rho,\eta)]^{-1}\xi=-{\partial_y  {\mathcal T}(\rho,\eta,\xi) \over
  \partial_{yy}\, \mathcal U(\rho,\eta) }\Big|_{J_\eta},
$$
  where $z=\mathcal T(\rho,\eta,\xi)$ is the solution of the problem
$$
   \mathcal A(\rho) z=z \quad  \mbox{in}\;\;
   \Bbb E_\eta,
   \qquad\;
   z =0 \quad  \mbox{on}\;\; \Gamma_s,
   \qquad\;
   z=\xi
  \quad  \mbox{on}\;\; J_\eta.
$$

   Notice that $D_\eta F(0,0)$ is  an isomorphism from $h^{m+\alpha}(\Bbb S)$ 
   onto $h^{m+1+\alpha}(\Bbb S)$ for all $m\in\Bbb N$, so classical 
   implicit function theorem in Banach spaces is not available here.  
   But on the other hand,
   similarly as (\ref{3.15}), we can show the mapping 
$$
    (\rho,\eta,\xi)\mapsto [D_\eta F(\rho,\eta)]^{-1}\xi  \; \mbox{  is smooth tame from }
   \; \mathcal O_\delta\times
   \widetilde{\mathcal O}_\delta^\infty\times C^\infty(\Bbb S) \; \mbox{ to } \;
    C^\infty (\Bbb S).
$$  
   Thus by Nash-Moser implicit function theorem (see Theorem 3.3.1 in Part III of
   [\ref{ham-82}]), there exist sufficiently small 
  $\delta_1, \delta'_1\in(0,r_0)$,
  and a unique smooth tame mapping $\mathcal S$ from  $\mathcal O_{\delta_1}$ 
  to  $\widetilde{\mathcal O}_{\delta'_1}^\infty$ 
   such that 
 $$
   \mathcal S(0)=0\qquad
   \mbox{and} \qquad
   F(\rho,\mathcal S(\rho))=0.
$$
  By letting $u=\mathcal U(\rho,\mathcal S(\rho))$ and $\eta=\mathcal S(\rho)$, 
  we see that $(u,\eta)$ is the solution of problem (\ref{3.9}), and the mapping 
  $\rho\mapsto(u,\eta)$ is smooth.
  The proof is complete. $\qquad\Box$
  
\medskip

\medskip

  By the proof of Lemma 3.2, for any $\rho\in \mathcal O_{\delta_1}$,
  problem (\ref{3.8}) has a unique solution 
\begin{equation}\label{3.17}
  u=\left\{
  \begin{array}{ll}
  \mathcal U(\rho,\mathcal S(\rho))\quad\;& \mbox{in}
  \;\; \Bbb E_{\mathcal S(\rho)},
  \\  [0.2 cm]
  \hat\sigma\qquad& \mbox{in}\;\;
   \Bbb D_{\mathcal S(\rho)},
  \end{array}
  \right.
  \qquad\quad\mbox{and}\qquad\quad
  \eta=\mathcal S(\rho).
\end{equation}
 Next we consider the following problem 
\begin{equation}\label{3.18}
\left\{
\begin{array}{rll}
  \mathcal A(\rho)v\,&=-\mu(u-\tilde{\sigma})\chi_{\Bbb E_\eta} 
  +\nu \chi_{\Bbb D_\eta}
  \quad  \qquad \; &\hbox{in}\;\;
   \Omega_s, 
\\ [0.2 cm]
  v\,&=\gamma\mathcal K(\rho)
  \qquad & \hbox{on}\;\; \Gamma_s,
\\ [0.2 cm]
   [\![v]\!]\,&=0,\quad \,[\![\partial_n v]\!]=0\quad
   & \hbox{on}\;\; J_{\eta}, 
\\ [0.2 cm]
  \partial_y v\,&=0 & \hbox{on}\;\; \Gamma_0,
\end{array}
\right.
\end{equation}
  where $u$ and $\eta$ are given by (\ref{3.17}).
  For the sake of simplicity,  
  we first study 
\begin{equation}\label{3.19}
\quad\;\,
\left\{
\begin{array}{rll}
  \mathcal A(\rho)w^+\,&=-\mu(\mathcal U(\rho,\mathcal S(\rho))-\tilde\sigma)
  \qquad\quad\;\; & \mbox{in} \;\; \Bbb E_{\mathcal S(\rho)},
  \\ [0.2 cm] 
  \mathcal A(\rho) w^-\,& = \nu   \qquad  &\mbox{in}\;\;\Bbb D_{\mathcal S(\rho)}, 
\\ [0.2 cm]
  w^+\,&=0
  \qquad & \hbox{on}\;\; \Gamma_s,
\\ [0.2 cm]
   w^+\,&=w^- 
   & \hbox{on}\;\; J_{\mathcal S(\rho)}, 
\\ [0.2 cm]
  \partial_n w^+\,&=\partial_n w^-\quad
   & \hbox{on}\;\; J_{\mathcal S(\rho)}, 
\\ [0.2 cm]
  \partial_y w^-\,&=0 & \hbox{on}\;\; \Gamma_0.
\end{array}
\right.
\end{equation}

\medskip

  {\bf Lemma 3.4} \ \  {\em There exists a positive constant $\delta_2\in (0,\delta_1)$
  such that for any $\rho\in \mathcal O_{\delta_2}$, problem $(\ref{3.19})$ has 
  a unique solution $(w^+,w^-)\in h^{4+\alpha}(\Bbb E_{\mathcal S(\rho)})\times
  h^{4+\alpha}(\Bbb D_{\mathcal S(\rho)})$, and the mapping 
  $\rho\mapsto (w^+,w^-)$ is smooth in $\mathcal O_{\delta_2}$.
  }
  
\medskip

  {\bf Proof.} \ \  For given $\rho\in\mathcal O_{\delta_1}$ and $\zeta\in h^{4+\alpha}(\Bbb S)$,
  we consider 
\begin{equation}\label{3.20}
\left\{
\begin{array}{rll}
  \mathcal A(\rho)w^+\,& =-\mu(\mathcal U(\rho,\mathcal S(\rho))-\tilde\sigma)
  \quad & \mbox{in}\;\; \Bbb E_{\mathcal S(\rho)},
  \\ [0.2 cm]
  w^+\,& =\zeta \qquad
  & \mbox{on}\;\; J_{\mathcal S(\rho)},
    \\ [0.2 cm]
    w^+\,& =0 \qquad& \mbox{on}\;\; \Gamma_s,
\end{array}
\right.
  \qquad\quad
\left\{
\begin{array}{rll}
  \mathcal A(\rho)w^-\,& =\nu
  \quad & \mbox{in}\;\; \Bbb D_{\mathcal S(\rho)},
  \\ [0.2 cm]
  w^-\,& =\zeta \quad& \mbox{on}\;\; J_{\mathcal S(\rho)},
  \\ [0.2 cm]
  \partial_y w^-\,& =0 \quad
  & \mbox{on}\;\; \Gamma_0.
\end{array}
\right.
\end{equation}
  From Lemma 3.2, we see $\mathcal S(\rho)\in C^{\infty}(\Bbb S)$
  and $\mathcal U(\rho,\mathcal S(\rho))\in h^{4+\alpha}(\Bbb E_{\mathcal S(\rho)})$. 
  By classical regularity theory of elliptic differential equations, 
  problem (\ref{3.20}) has a unique 
  solution $(w^+,w^-)$ such that
\begin{equation}\label{3.21}
  w^+:=\mathcal W^+(\rho,\zeta)
  \in h^{4+\alpha}(\Bbb E_{\mathcal S(\rho)})\qquad
  \mbox{ and }\qquad
  w^-:=\mathcal W^-(\rho,\zeta)\in h^{4+\alpha}(\Bbb D_{\mathcal S(\rho)}).
\end{equation}
  Since the mappings $\mathcal S$ and $\mathcal U$ are both
  smooth in $\mathcal O_{\delta_1}$, the mappings $\mathcal W^+$
  and $\mathcal W^-$
  are also smooth in $\mathcal O_{\delta_1}\times h^{4+\alpha}(\Bbb S)$. 

  Recall that $\partial_n$ is outward normal derivative on $J_{\mathcal S(\rho)}$
  with respect to $\Bbb E_{\mathcal S(\rho)}$.
  Define a mapping
  $G: \mathcal O_{\delta_1}\times h^{4+\alpha}(\Bbb S)\to h^{3+\alpha}(\Bbb S)$ by
\begin{equation}\label{3.22}
  G(\rho,\zeta)=
  \partial_n \mathcal W^+(\rho,\zeta)\Big |_{J_{\mathcal S(\rho)}}
  -\partial_n \mathcal W^-(\rho,\zeta)\Big |_{J_{\mathcal S(\rho)}}\qquad
  \mbox{for}\;\;\rho\in\mathcal O_{\delta_1},\;\;
  \zeta\in h^{4+\alpha}(\Bbb S).
\end{equation}
  It is easy to see that problem (\ref{3.19}) is equivalent to the equation
  $G(\rho,\zeta)=0$.  
  
  Since $\mathcal W^+$ and $\mathcal W^-$ are smooth, 
  we have
\begin{equation}\label{3.23}
  G\in C^\infty\big(\mathcal O_{\delta_1}\times h^{4+\alpha}(\Bbb S),
  h^{3+\alpha}(\Bbb S)\big).
\end{equation}
  By (\ref{2.1})--(\ref{2.3}), we see $G(0, p_0)=0$, where $p_0=p_s(\eta_s)$.  
  Note that 
$$
  \mathcal S(0)=0, \quad
  \mathcal U(0,\mathcal S(0))=\sigma_s,
  \quad
  \mathcal W^+(0,p_0)=p_s\big|_{\Bbb E_s}, 
  \quad 
 \mathcal W^-(0,p_0)=p_s\big|_{\Bbb D_s}.
$$ 
  By a direct computation,
  we have 
$$
  D_\zeta G(0,p_0)\xi =- \partial_y z^+\big|_{J_s}+\partial_y z^-\big|_{J_s}
  \qquad
  \mbox{for}\quad \xi\in h^{4+\alpha}(\Bbb S),
$$
  where $z^+$ and $z^-$ are the solutions to the following two problems, respectively,
\begin{equation}\label{3.24}
\begin{array}{l}
  \Delta z^+ =0
  \quad \mbox{in}\;\; \Bbb E_s,
  \qquad
  z^+ =\xi \quad \mbox{on}\;\; J_s,
  \qquad    
  z^+ =0 \quad \mbox{on}\;\; \Gamma_s,
\\ [0.3 cm]
  \Delta z^- =0  \quad  \mbox{in}\;\; \Bbb D_s,
  \qquad
  z^-=\xi \quad \mbox{on}\;\; J_s,
  \qquad
  \partial_y z^- =0 \quad\mbox{on}\;\; \Gamma_0.
\end{array}
\end{equation}
  For any $\xi\in C^\infty(\Bbb S)$ with the expression 
  $\displaystyle \xi(x)=\sum_{k\in \Bbb Z} \xi_k e^{ikx}$,
  we obtain
\begin{equation}\label{3.25}
  D_\zeta G(0,p_0)\xi = \sum_{k\in \Bbb Z} \tau_k \xi_k  e^{ikx},
\end{equation}
  where $\tau_0=(\rho_s-\eta_s)^{-1}$ and
  $\tau_k=k(\coth k(\rho_s-\eta_s)+\tanh k\eta_s )$ for
  $k\neq 0$, $k\in\Bbb Z$.
  
  Obviously, there exist two positive constant $C_1$ and $C_2$ such that
$$
   C_1\sqrt{k^2+1}\le \tau_k\le C_2 \sqrt{k^2+1}.
$$
  It implies that
\begin{equation}\label{3.26}
  D_\zeta G(0,p_0) \mbox{ is an isomorphism from }
  H^{r+1}(\Bbb S) \mbox{  onto } H^r(\Bbb S)
  \mbox{  for }r>0,
\end{equation}
  where
  $\displaystyle H^r(\Bbb S)=\{f\in L^2(\Bbb S): 
  \sum_{k\in\Bbb Z} (k^2+1)^r |\widehat{f}(k)|^2< +\infty\}$.
  
  From (\ref{3.25}), we easily obtain that for $\xi\in C^\infty(\Bbb S)$ 
  with  $\displaystyle \xi(x)=\sum_{k\in \Bbb Z} \xi_k e^{ikx}$,
$$
  [D_\zeta G(0,p_0)]^{-1} \xi = \sum_{k\in \Bbb Z} \tau_k^{-1} \xi_k  e^{ikx}.
$$
  Define a function $\tau(x)=x(\coth (\rho_s-\eta_s)x+\tanh \eta_sx )$
  for $|x|\ge1$. It is easy to verify that
$$
  \tau_k=\tau(k)\quad \mbox {for }\; k\neq 0\quad \mbox{and}\quad
    \sup_{|x|\ge1} |\tau'(x)|+|x||\tau''(x)|<+ \infty.
$$
  Using above relations one can prove that
$$
\left\{
\begin{array}{l}
  \displaystyle \sup_{k\in \Bbb Z} | k |\big|{1\over \tau_k}\big| < +\infty,
  \\ [0.4 cm]
   \displaystyle  \sup_{k\in \Bbb Z} | k |^2\big|{1\over \tau_{k+1}}
   -{1\over \tau_k} \big| < +\infty,
  \\ [0.4 cm]
   \displaystyle  \sup_{k\in \Bbb Z} | k |^3\big|
   {1\over \tau_{k+2}}-{2\over \tau_{k+1}}+
   {1\over \tau_k}\big| < +\infty.
\end{array}
\right.
$$
  Then by Theorem 4.5 of [\ref{are-bu-04}] (or  [\ref{sch-tri-87}]),
  we have 
\begin{equation}\label{3.27}
  [D_\zeta G(0,p_0)]^{-1}\in L(C^r(\Bbb S), C^{r+1}(\Bbb S))\qquad
  \mbox{for }\;\; r>0.
\end{equation}
  By Sobolev embedding theorem, $ H^{4+r}(\mathbb{S})\,
   {\hookrightarrow}\, C^{4+\alpha}(\mathbb{S})$ for 
   $r>{3/2}$.
  Notice that $h^{4+\alpha}(\mathbb{S})$ is the closure of 
  $H^{4+r}(\mathbb{S})$  in $C^{4+\alpha} (\mathbb{S})$ for 
  $r>{3/2}$.
  By  (\ref{3.26}) and (\ref{3.27}), we obtain that
  $[D_\zeta G(0,p_0)]^{-1}\in L(h^{3+\alpha}(\Bbb S), h^{4+\alpha}(\Bbb S))$
  and 
$$
  D_\zeta G(0,p_0) \mbox{ is an isomorphism from }
  h^{4+\alpha}(\Bbb S) \mbox{  onto } h^{3+\alpha}(\Bbb S).
$$
  Hence by classical implicit function theorem
  in Banach spaces,  there exist sufficiently small 
  constants $\delta_2, \delta_2'\in(0,\delta_1)$,
  and a unique mapping $\mathcal R\in C^\infty(
  \mathcal O_{\delta_2},  h^{4+\alpha}(\Bbb S))$ 
  such that 
 $$
   \mathcal R(0)=p_0,\qquad
   \|\mathcal R(\rho)-p_0\|_{h^{4+\alpha}(\Bbb S)}\le \delta'_2 
   \qquad \mbox{and} \qquad
   G(\rho,\mathcal R(\rho))=0.
$$
  By letting $(w^+,w^-)=(\mathcal W^+(\rho,\mathcal R(\rho)), 
  \mathcal W^-(\rho,\mathcal R(\rho)))$, we see that $(w^+,w^-)$
   is the solution of problem
  (\ref{3.19}),  and the desired result follows immediately.  
  $\qquad \Box$

\medskip

\medskip

\medskip
  
  By the proof of Lemma 3.4,  we denote
\begin{equation}\label{3.28}
  \mathcal W(\rho)=\left\{
  \begin{array}{l}
  \mathcal W^+(\rho,\mathcal R(\rho))\qquad\mbox{in}\;\; \Bbb E_{\mathcal S(\rho)},
  \\ [0.2 cm] 
  \mathcal W^-(\rho,\mathcal R(\rho))\qquad\mbox{in}\;\; \Bbb D_{\mathcal S(\rho)},
  \end{array}
  \right.
  \qquad\;\;\mbox{for}\;\;\;\; \rho\in\mathcal O_{\delta_2}.
\end{equation}
  Consider the problem
\begin{equation}\label{3.29}
  \mathcal A(\rho) v_0=0
   \quad  \hbox{in}\;\;
   \Omega_s, 
\qquad
  v_0 =\gamma\mathcal K(\rho)
  \quad  \hbox{on}\;\; \Gamma_s,
\qquad
  \partial_y v_0 =0\quad  \hbox{on}\;\; \Gamma_0.
\end{equation}
  Note that by (\ref{3.3}), we have
\begin{equation}\label{3.30}
  \mathcal K \in C^\infty(\mathcal O_{\delta_2},h^{2+\alpha}(\Bbb S)).
\end{equation}
   By classical regularity theory of elliptic differential equations,
  problem (\ref{3.29}) has a unique solution
  $v_0:= \mathcal V(\rho)\in h^{2+\alpha}(\Omega_s)$. Moreover,
  by (\ref{3.5}), (\ref{3.30}) and Lemma 2.3 in [\ref{esc-sim-97-1}], 
\begin{equation}\label{3.31}
  \mathcal V \in C^\infty (\mathcal O_{\delta_2}, h^{2+\alpha}(\Omega_s)).
\end{equation}
  From (\ref{3.28}), (\ref{3.29}) and Lemma 3.4, for any $\rho\in\mathcal O_{\delta_2}$,
   we  see problem (\ref{3.18})
  has a unique solution 
\begin{equation}\label{3.32}
  v= \mathcal V(\rho)+ \mathcal W(\rho).
\end{equation}
  
  Later on, we always fix $0<\delta\le\delta_2$.
  Define a mapping $\Psi: \mathcal O_\delta\to h^{1+\alpha}(\Bbb S)$ by
\begin{equation}\label{3.33}
  \Psi(\rho):=\mathcal B(\rho)\mathcal V(\rho) +\mathcal B(\rho)\mathcal W(\rho)
  \qquad\mbox{for}\;\; \rho\in\mathcal O_\delta.
\end{equation}
  It follows from $(\ref{3.5})$, (\ref{3.31}) and Lemma 3.4 that
\begin{equation}\label{3.34}
  \Psi  \in C^\infty(\mathcal O_\delta, h^{1+\alpha}(\Bbb S)).
\end{equation}
  With all above reductions, we see that problem (\ref{3.7}) is equivalent to
  the following Cauchy problem
\begin{equation}\label{3.35}
\left\{
\begin{array}{ll}
  \partial_t \rho  + \Psi(\rho)=0
  \qquad\;\; & \mbox{on}\;\; \Bbb S, \;\; t>0,
  \\ [0.2 cm]
  \rho(0)=\rho_0   \qquad & \mbox{on}\;\; \Bbb S.
\end{array}
\right.
\end{equation}

   More precisely, we have

 \medskip

{\bf Lemma 3.5} \ \ {\em  The function $\rho$ is the
  solution of problem $(\ref{3.35})$ if and only if 
  $(u, v, \eta, \rho)$ is the solution of problem $(\ref{3.7})$
  with $(u, v,\eta)$ given by $(\ref{3.17})$ and $(\ref{3.32})$.
}

\medskip

  Next we study local well-posedness of problem (\ref{3.35}).
  For any $\rho\in\mathcal O_\delta$, we define the Fr\'echet 
  derivative of nonlinear operator $\Psi$ at $\rho$ by 
$$
  \displaystyle D \Psi(\rho)\zeta:=\lim_{\varepsilon\to0}
  {\Psi(\rho+\varepsilon\zeta)-\Psi(\rho)\over\varepsilon}
  \qquad\quad \mbox{for}\;\; \zeta\in h^{4+\alpha}(\Bbb S).
$$
  
    Let  $E_0$ and $E_1$ be two Banach spaces,   $E_1$
    is densely and continuously embedded into $E_0$.  Denote by 
    $\mathcal H(E_1,E_0)$  the subspace of all linear operators
    $A\in L(E_1,E_0)$ such that $-A$ generates a strongly continuous
    analytic semigroup on $E_0$.  We have the following result:

\medskip

{\bf Lemma 3.6} \ \  {\em 
  $D \Psi (\rho)\in \mathcal H(h^{4+\alpha}(\Bbb S),h^{1+\alpha}(\Bbb S))$
  for $\rho\in\mathcal O_\delta$.  }
\medskip

  {\bf Proof.} \ \ Let $\Psi_1(\rho):=\mathcal B(\rho)\mathcal V(\rho)$
  and $\Psi_2:=\mathcal B(\rho)\mathcal W(\rho)$,   then
$$
  \Psi(\rho)=\Psi_1(\rho)+\Psi_2(\rho)\qquad \mbox{for}\;\; \rho\in\mathcal O_\delta.
$$ 
  Notice that the following problem is the corresponding transformed 
  periodic Hele-Shaw model with surface tension:
$$
 \left\{
\begin{array}{rll}
  \mathcal A(\rho) v_0\,&=0
  \quad  \qquad  &\hbox{ in }\;
   \Omega_s, \;t>0,
\\ [0.2 cm]
  \partial_y v_0\,&=0 & \hbox{ on }\; \Gamma_0,\;t>0,
\\ [0.2 cm]
  v_0 \,&=\gamma\mathcal K(\rho)
  \qquad & \hbox{ on }\; \Gamma_s,\; t>0,
\\ [0.2 cm]
   \rho_t\,&=-\mathcal B(\rho) v_0 \qquad \qquad 
   \qquad & \hbox{ on }\; \Bbb S, \; \; \, t>0,
\end{array}
\right.
$$
  and similarly,  it can be reduced to $\partial_t\rho+\Psi_1(\rho)=0$ for $t>0$.
  Thus by well-known results of Hele-Shaw models (cf. [\ref{esc-sim-97-1}]),
  we have $D\Psi_1(\rho)\in \mathcal H(h^{4+\alpha}(\Bbb S),h^{1+\alpha}(\Bbb S))$,
  for any $\rho\in\mathcal O_\delta$.
  
  On the other hand, by Lemma 3.4, (\ref{3.5}) and (\ref{3.28}), we have 
  $\Psi_2\in C^\infty(\mathcal O_\delta,h^{3+\alpha}(\Bbb S))$ and 
  $D\Psi_2(\rho)\in L(h^{4+\alpha}(\Bbb S),h^{3+\alpha}(\Bbb S))$. Since
  $h^{3+\alpha}(\Bbb S)$ is compactly embedded into $h^{1+\alpha}(\Bbb S)$,
  by the well-known perturbation result (cf. Theorem I.1.5.1 in [1], or 
  Proposition 2.4.3 in [\ref{lunardi}]), we get the desired result. $\qquad\Box$

\medskip
  
  The above result implies that problem (\ref{3.35}) is of parabolic type
  in $\mathcal O_\delta$.
  Thus by using analytic semigroup theory and applications to parabolic
  differential problems (see [\ref{amann}] and [\ref{lunardi}]), we get the local
  well-posedness.

\medskip

{\bf Theorem 3.7} \ \  {\em Given $\rho_0\in \mathcal O_\delta$. There exists
  a maximal $T>0$ such that problem $(\ref{3.35})$ has a unique solution
  $\rho\in C([0,T),\mathcal O_\delta)\cap C^1([0,T), h^{1+\alpha}(\Bbb S))$.
 }

\medskip

  From Theorem 3.7, and combining Lemma 3.1, Lemma 3.5, we see that free
  boundary problem (\ref{1.1}) is locally wellposed, and given $\rho_0\in \mathcal O_\delta$,
  there exists a unique solution $(\sigma, p,\eta, \rho)$ of problem (\ref{1.1}).

\medskip

\hskip 1em

\section{Linearization and Eigenvalues}
\setcounter{equation}{0}
\hskip 1em

  In this section we study linearization of problem (\ref{3.35})
  at the stationary solution $\rho=0$, and compute all eigenvalues
  of $D\Psi(0)$.

  First, we study the linearization
  of free boundary problem (\ref{1.1}) at flat stationary 
  solution $(\sigma_s, p_s, \eta_s, \rho_s)$.
   Let
\begin{equation}\label{4.1}
  \sigma=\sigma_s+\epsilon \phi (x,y,t),\quad
  p=p_s+\epsilon \psi(x,y,t),\quad
  \eta=\eta_s+\epsilon \xi(x,t), \quad
  \rho=\rho_s+\epsilon \zeta(x,t),
\end{equation}
  where $\phi$, $\psi$, $\xi$ and $\zeta$ are unknown functions. 
  At each time $t>0$,
   by (\ref{3.6}), the mean curvature of the curve
   $y=\rho_s+\epsilon\zeta$ can be expressed  by
\begin{equation}\label{4.2}
  \mathcal K(\epsilon\zeta)=-\epsilon\zeta_{xx} +O(\epsilon^2).
\end{equation}
  Let ${\bf n}_{\epsilon\zeta}=(-\epsilon\zeta_x,1)$ be the
  outward normal  
  direction on $y=\rho_s+\epsilon\zeta$. 
  We compute 
\begin{equation}\label{4.3}
\begin{array}{rl}
  \langle\nabla p, {\bf n}_{\epsilon\zeta}\rangle \big|_{y=\rho_s+\epsilon\zeta}
   &=(-\epsilon p_x \zeta_x + p_y )\big|_{y=\rho_s+\epsilon\zeta}
  = \partial_y(p_s+\epsilon\psi)\big|_{y=\rho_s+\epsilon\zeta}
  +O(\epsilon^2)
  \\ [0.3 cm]
  &=\epsilon p_s''(\rho_s)\zeta+\epsilon \partial_y\psi\big|_{y=\rho_s}+O(\epsilon^2)
  \\ [0.3 cm]
  &  =-\epsilon\big[\mu(\bar\sigma-\tilde\sigma)\zeta-\partial_y\psi\big|_{y=\rho_s}
  \big]+O(\epsilon^2).
  \end{array}
\end{equation}
  By substituting (\ref{4.1}) into problem (\ref{1.1}),  collecting all 
  first order $\epsilon$-terms and with the aid of  (\ref{4.2}), (\ref{4.3})
  and the fact that 
$$
  \sigma_s''(\eta_s^+)=\hat\sigma, 
  \qquad 
  \sigma_s''(\eta_s^-)=0,
  \qquad 
  \sigma_s'(\rho_s)=\sqrt{\bar\sigma^2-\hat\sigma^2},
$$
$$
  p_s''(\eta_s^+)=-\mu(\hat\sigma-\tilde\sigma),
  \qquad
  p_s''(\eta_s^-)=\nu,
  \qquad
  p_s'(\rho_s)=0,
$$ 
  we obtain the linearization of problem (\ref{1.1}) at  
  $(\sigma_s, p_s, \eta_s, \rho_s)$
  is given by
\begin{equation}\label{4.4}
\left\{
\begin{array}{ll}
   \Delta \phi=\phi \chi_{ \Bbb E_s}
 \qquad & \mbox{in}\;\;  \Omega_s,\;t>0,
\\ [0.3 cm]
   \Delta \psi=-\mu\phi \chi_{\Bbb E_s}
 \qquad & \mbox{in}\;\;  \Omega_s,\;t>0,
\\ [0.3 cm]
  \phi=-\sqrt{\bar\sigma^2-\hat\sigma^2}\zeta,\qquad 
  \psi =-\gamma\zeta_{xx}
  \qquad\qquad & \mbox{on}\;\;  \Gamma_s,\;t>0,
\\ [0.3 cm]
   \phi=0,\qquad  [\![\partial_y \phi]\!]=-\hat\sigma\xi
   \qquad &  \mbox{on}\;\;  J_s,\;t>0,
\\ [0.3 cm]
   [\![\psi]\!]=0,\quad\;  [\![\partial_y \psi]\!]=
   \big(\mu(\hat\sigma-\tilde\sigma)+\nu\big)\xi
   \qquad &  \mbox{on}\;\;  J_s,\;t>0,
\\ [0.3 cm]
   \partial_y\phi=0 ,\quad \partial_y\psi=0  & \mbox{on}\;\;  \Gamma_0,\;t>0,
\\ [0.3 cm]
   \partial_t\zeta=-\partial_y \psi\big|_{y=\rho_s} +\mu(\bar\sigma-\tilde\sigma)\zeta
   & \mbox{on}\;\;  \Bbb S,\;\,\;t>0.
\end{array}
\right.
\end{equation}
   For any given $\zeta\in h^{4+\alpha}(\Bbb S)$, by solving problem 
   $(4.4)_1$--$(4.4)_6$, we get a unique solution $(\phi,\psi,\xi)$.
   
   Since problem (\ref{1.1}) is equivalent to problem (\ref{3.35}), their corresponding
   linearizations at flat stationary solution are also equivalent. It implies that 
\begin{equation}\label{4.5}
  D\Psi(0)\zeta = \partial_y \psi\big|_{y=\rho_s} - \mu(\bar\sigma-\tilde\sigma)\zeta
  \qquad\mbox{for}\;\; \zeta\in h^{4+\alpha}(\Bbb S).
\end{equation}
   
   Next, we given an explicit expression of $D\Psi(0)$ and study its eigenvalues.
   For any given
\begin{equation}\label{4.6}
  \zeta(x)=\displaystyle \sum_{k\in\Bbb Z} c_k {\bf e}^{ikx}
   \in C^\infty(\Bbb S),
\end{equation}
  set
 \begin{equation}\label{4.7}
  \phi(x,y)=\sum_{k\in\Bbb Z} a_k(y) {\bf e}^{ikx},
  \qquad
  \psi(x,y)=\sum_{k\in\Bbb Z} b_k(y) {\bf e}^{ikx},
  \qquad
  \xi(x)=\displaystyle \sum_{k\in\Bbb Z} d_k {\bf e}^{ikx},
\end{equation}
  where $a_k(y)$ and $b_k(y)$ are unknown functions, $d_k$ is unknown
  coefficient for each $k\in\Bbb Z$.
  
  Substituting (\ref{4.6}) and (\ref{4.7}) into (\ref{4.4}), we see that 
  for each $k\in\Bbb Z$,  there hold 
\begin{equation}\label{4.8}
  \left\{
  \begin{array}{ll}
  \displaystyle a_{k}''-k^2 a_{k}=a_{k} \qquad
  \qquad \qquad\,  & \mbox{for}\;\; \eta_s< y <\rho_s,
\\ [0.4 cm]
  a_k(y)=0 \qquad
  & \mbox{for}\;\; 0< y \le \eta_s,
\\ [0.4 cm]
  \displaystyle  a_k' (\eta_s^+)=-\hat\sigma d_k,
\\ [0.4 cm]
  a_k(\rho_s)=-\sqrt{\bar\sigma^2-\hat\sigma^2}c_k,
 \end{array}
 \right.
\end{equation}
  and
\begin{equation}\label{4.9}
\left\{
\begin{array}{l}
   \displaystyle  b_{k}''-k^2 b_{k} =- \mu a_{k}
   \qquad\qquad\quad\;\;\;\;
   \mbox{for}\;\; \eta_s< y <\rho_s,
\\ [0.4 cm]
  \displaystyle  b_{k}''-k^2 b_{k} = 0
   \qquad\qquad\qquad\qquad
   \mbox{for}\;\; 0< y <\eta_s,
\\ [0.4 cm]
  \displaystyle   b_k'(\eta_s^+)
   = b_k'(\eta_s^-)+\big(\mu(\hat\sigma-\tilde\sigma)+\nu\big) d_k,
\\ [0.4 cm]
  b_k(\eta_s^+)=b_k(\eta_s^-),
\\ [0.4 cm]
  b_k (\rho_s)=\gamma k^2 c_k,
\\ [0.4 cm]
  \displaystyle  b_k'(0)=0.
\end{array}
\right. \;
\end{equation}
   By solving problem (\ref{4.8}), we obtain that for each $k\in\Bbb Z$, 
\begin{equation}\label{4.10}
 a_k(y)=\left\{
 \begin{array}{ll}
 \displaystyle
 -{\sinh \sqrt{k^2+1}(y-\eta_s)\over
 \sinh \sqrt{k^2+1}(\rho_s-\eta_s)}\sqrt{\bar\sigma^2-\hat\sigma^2}c_k \qquad
 & \mbox{for}\;\; \eta_s\le y\le \rho_s,
 \\[ 0.4 cm]
 0 \qquad\quad  & \mbox{for}\;\; 0< y< \eta_s,
 \end{array}
 \right.
\end{equation}
  and
\begin{equation}\label{4.11}
   d_k=
   { \sqrt{k^2+1}\sqrt{\bar\sigma^2-\hat\sigma^2}c_k 
   \over \hat\sigma
   \sinh \sqrt{k^2+1}(\rho_s-\eta_s)}.
\end{equation}
  Then by solving problem (\ref{4.9}), an elementary computation shows that,
  for  each $k\neq0$, $k\in\Bbb Z$,
\begin{equation}\label{4.12} 
  b_k(y)=\left\{
  \begin{array}{ll}
  \displaystyle 
  -\mu a_k(y)+(\gamma k^2-\mu\sqrt{\bar\sigma^2-\hat\sigma^2})c_k
  {\cosh ky\over \cosh k\rho_s}+ e_k {\sinh k(\rho_s-y)\over
  \sinh k(\rho_s-\eta_s)}\quad& \mbox{for}\;\; \eta_s\le y\le \rho_s,
  \\ [0.6 cm]
  \displaystyle 
  (\gamma k^2-\mu\sqrt{\bar\sigma^2-\hat\sigma^2})c_k
  {\cosh ky\over \cosh k\rho_s}+ e_k {\cosh ky\over
  \cos k\eta_s}\qquad& \mbox{for}\;\; 0<y<\eta_s,
  \end{array}
  \right.
\end{equation}
  where
\begin{equation}\label{4.13}
  e_k={(\mu\tilde\sigma-\nu) d_k\over
  k[\coth k(\rho_s-\eta_s)+\tanh k\eta_s]}\qquad
  \qquad \mbox{for}\;\ k\neq0,\; k\in\Bbb Z.
\end{equation}
  By using (\ref{2.4}) and (\ref{4.11}),  we have $d_0=c_0$, then
\begin{equation}\label{4.14}
  b_0(y)=\left\{
  \begin{array}{ll}
  \displaystyle 
  \Big[\mu\sqrt{\bar\sigma^2-\hat\sigma^2}
  \Big({\sinh(y-\eta_s)\over \sinh(\rho_s-\eta_s)}-1\Big)
  +(\mu\tilde\sigma-\nu)
  (\rho_s-y)\Big]c_0 
  \quad\;\; \; & \mbox{for}\;\; \eta_s\le y\le \rho_s,
  \\ [0.5 cm]
  \displaystyle 
  \Big[ -\mu\sqrt{\bar\sigma^2-\hat\sigma^2}+(\mu\tilde\sigma-\nu)
  (\rho_s-\eta_s)\Big]c_0 \qquad& \mbox{for}\;\; 0<y<\eta_s.
  \end{array}
  \right.
\end{equation}
   By (\ref{4.10})--(\ref{4.13}), for $k\neq 0$, we compute
\begin{equation}\label{4.15}
  \begin{array}{rl}
     &b_k'(\rho_s)-\mu (\bar\sigma-\tilde\sigma)c_k  
     \\ [0.3 cm]
  = \,&  \displaystyle
  -\mu a_k'(\rho_s)+(\gamma k^2-\mu
  \sqrt{\bar\sigma^2-\hat\sigma^2})c_k k\tanh k\rho_s
  -  {k\, e_k\over\sinh k(\rho_s-\eta_s)}-\mu (\bar\sigma-\tilde\sigma)c_k
  \\ [0.3 cm]
  = \, & \lambda_k(\gamma) c_k,
  \end{array}
\end{equation}
  where   
\begin{equation}\label{4.16}  
  \begin{array}{rl}
  \lambda_k(\gamma) 
  \,& = \gamma k^3 \tanh k\rho_s + \mu\sqrt{\bar\sigma^2-\hat\sigma^2}\big[
  \sqrt{k^2+1}\coth \sqrt{k^2+1}(\rho_s-\eta_s)- k\tanh k\rho_s\big]
  \\ [0.3 cm]
   & \displaystyle 
   +{(-\mu\tilde\sigma+\nu)\sqrt{k^2+1}  \sqrt{\bar\sigma^2-\hat\sigma^2}  \over
   \hat\sigma\sinh k(\rho_s-\eta_s)\sinh \sqrt{k^2+1} (\rho_s-\eta_s)
   \big[ \coth k(\rho_s-\eta_s)+\tanh k\eta_s\big]}-\mu (\bar\sigma-\tilde\sigma),
  \end{array}  
\end{equation}
  for $k\neq 0$ and $\gamma>0$.

  Note that (\ref{2.4}) implies  $\coth (\rho_s-\eta_s)=\bar\sigma/ 
  \sqrt{\bar\sigma^2-\hat\sigma^2}$. 
  Then from (\ref{4.14}) we compute
\begin{equation}\label{4.17}
  \begin{array}{rl}
   b_0'(\rho_s)-\mu(\bar\sigma-\tilde\sigma)c_0
  \, & =  \mu c_0 \sqrt{\bar\sigma^2-\hat\sigma^2}
  \coth (\rho_s-\eta_s) -(\mu\tilde\sigma- \nu)c_0 
  -\mu(\bar\sigma-\tilde\sigma)c_0
  \\ [0.3 cm]
  \, &=  \mu \bar\sigma c_0- (\mu\tilde\sigma-\nu)c_0
  -\mu(\bar\sigma-\tilde\sigma)c_0
  \\ [0.3 cm]
  \, & = \nu c_0.
 \end{array}
\end{equation}

   By (\ref{4.5})--(\ref{4.7}) and (\ref{4.15})--(\ref{4.17}), 
   we have 

\medskip

\medskip

  {\bf Lemma 4.1} \ \  {\em For any $\zeta\in
  C^\infty(\Bbb S)$ given by $\displaystyle\zeta=\sum_{k\in\Bbb Z}
  c_{k}{\bf e}^{ikx}$, there holds
\begin{equation}\label{4.18}
  D\Psi(0)\zeta=\sum_{k\in\Bbb Z}
  \lambda_k(\gamma)c_{k}{\bf e}^{ikx},
\end{equation}
  where $\lambda_k(\gamma)$ is given by $(\ref{4.16})$ for $k\neq0$,
  and $\lambda_0(\gamma) \equiv \nu$.
}

\medskip

\medskip

   Obviously, for each $k\in\Bbb Z$ and $\gamma>0$, 
   $\lambda_k(\gamma)$ is an eigenvalue of the linearized 
   operator $D\Psi(0)$.  We have the following properties:
 \medskip
 
   {\bf Lemma 4.2} \ \    {\em $(i)$ For any $\gamma>0$,
   $\lim_{k\to\infty}\lambda_k(\gamma)=+\infty$.
   
   $(ii)$ There exists a constant $\gamma_*>0$, such that
   if $\gamma>\gamma_*$, we have $\lambda_k(\gamma)>0$ for
   all $k\in \Bbb Z$; and if $0<\gamma<\gamma_*$, there exists
   at least an integer $k_0\in\Bbb Z$ such that $\lambda_{k_0}(\gamma)<0$. 
   }
   
\medskip
   
    {\bf Proof}. \ \ $(i)$ By a direct analysis, we have 
$$
  \lim_{k\to+\infty}\tanh k\rho_s=\lim_{k\to+\infty}\coth k(\rho_s-\eta_s)=1,
$$
$$
  \lim_{k\to-\infty}\tanh k\rho_s=\lim_{k\to-\infty}\coth k(\rho_s-\eta_s)=-1,
$$    
$$
  \lim_{k\to \infty} \big(\sqrt{k^2+1}\coth\sqrt{k^2+1}(\rho_s-\eta_s)- k\tanh k\rho_s
  \big)=0.
$$    
 Hence by (\ref{4.16}), we immediately obtain $\lim_{k\to\infty}\lambda(\gamma)
 =+\infty$ for any $\gamma>0$.
 
 $(ii)$ Define a sequence $\{\gamma_k\}_{k\neq0}$ by
\begin{equation}\label{4.19}  
    \begin{array}{rl}
  \gamma_k 
  \,&\displaystyle := {1\over k^3 \tanh k\rho_s }\Big\{ \mu\sqrt{\bar\sigma^2-\hat\sigma^2}
  \Big[k\tanh k\rho_s -\sqrt{k^2+1}\coth \sqrt{k^2+1}(\rho_s-\eta_s)\Big]
  \\ [0.3 cm]
   & \displaystyle 
   +{(\mu\tilde\sigma-\nu)\sqrt{k^2+1}  \sqrt{\bar\sigma^2-\hat\sigma^2}  \over
   \hat\sigma\sinh k(\rho_s-\eta_s)\sinh \sqrt{k^2+1} (\rho_s-\eta_s)
   \big[ \coth k(\rho_s-\eta_s)+\tanh k\eta_s\big]}+\mu (\bar\sigma-\tilde\sigma)
   \Big\}.
  \end{array}  
\end{equation}
  Clearly, we have
\begin{equation}\label{4.20}
 \lim_{k\to\infty}\gamma_k=0
 \qquad \mbox{and}\qquad
 \lim_{k\to\infty} k^3 \tanh k\rho_s \gamma_k =\mu(\bar\sigma-\tilde\sigma)>0.
\end{equation}
  
  Let 
\begin{equation}\label{4.21}
  \gamma_*:=\sup_{k\neq0}\{\gamma_k\}.  
\end{equation}
  By (\ref{4.20}) we see that $\gamma_*$ is well-defined and $\gamma_*>0$.
  
   By (\ref{4.19}),  we rewrite (\ref{4.16}) as 
\begin{equation}\label{4.22}
  \lambda_k(\gamma)=k^3  \tanh k\rho_s\, \big(\gamma-\gamma_k\big)\qquad
  \mbox{for}\;\; k\neq 0,\;k\in\Bbb Z.
\end{equation}
    Then the desired result follows from (\ref{4.20}) and (\ref{4.21}).
     $\qquad\Box$

 \medskip

   Denote $\sigma(D\Psi(0))$ by
   the spectrum of $D\Psi(0)$.
    By Lemma 4.1 and 
   Lemma 4.2, we have 

\medskip

{\bf Corollary 4.3} \ \ {\em $(i)$ If $\gamma>\gamma_*$, there exists
  a constant $\varpi>0$ such that
$$
  \sigma(D\Psi(0))\subset \{\lambda\in \Bbb C: {\rm Re} \lambda \ge \varpi\}.
$$
  
  $(ii)$ If $0<\gamma<\gamma_*$, then $\sigma(D\Psi(0))\cap
  \{\lambda\in \Bbb C: {\rm Re}\lambda<0\}\neq \emptyset$.

}

\medskip

  {\bf Proof}. \ \  Since $D\Psi(0)\in L(h^{4+\alpha}(\Bbb S),
   h^{1+\alpha}(\Bbb S))$, and  $h^{4+\alpha}(\Bbb S)$ is 
   compactly embedded into $h^{1+\alpha}(\Bbb S)$,
    we see that $\sigma(D\Psi(0))$
   consists of all eigenvalues. By Lemma 4.1, we easily show that
   all eigenvalues
   of the restriction of $D\Psi(0)$ in $H^{4+r}(\Bbb S)$ 
   are given by $\lambda_k(\gamma)$ for $k\in\Bbb Z$. 
   Since $h^{4+\alpha}(\Bbb S)$
   is the closure of $H^{4+r}(\Bbb S)$ in $C^{4+\alpha}(\Bbb S)$ for
   $r>3/2$, we have
$$
  \sigma(D\Psi(0))=\{\lambda_k(\gamma);k\in \Bbb Z\}.
$$
  
  Let $\gamma>\gamma_*$. By (\ref{4.21}) and (\ref{4.22}),
  we see that
$$
  \lambda_k(\gamma)\ge k^3\tanh k\rho_s \,(\gamma-\gamma_*)
  \ge \tanh \rho_s\, (\gamma-\gamma_*)>0 \qquad
  \mbox{for}\;\; k\neq 0.
$$
  Notice that $\lambda_0(\gamma)\equiv \nu>0$. Take 
  $\varpi\in (0,\min\{\tanh \rho_s\, (\gamma-\gamma_*), \nu\})$, then
  $\lambda(\gamma)>\varpi$ for all $k\in\Bbb Z$. It implies
  that the assertion $(i)$ holds.  The assertion $(ii)$ 
  directly follows from Lemma 4.2 $(ii)$. The proof is complete. 
   \qquad $\Box$

\medskip

\medskip
\hskip 1em

\section{Asymptotic stability}
\setcounter{equation}{0}
\hskip 1em

  In this section we study asymptotic stability of the 
  stationary solution $\rho=0$ of problem (\ref{3.35}) 
  and give a proof of our main result Theorem 1.2.

  Since problem (\ref{3.35}) is of parabolic type in 
  $h^{1+\alpha}(\Bbb S)$, by using geometric theory 
  of parabolic equations in Banach spaces, we have

\medskip

 {\bf Theorem 5.1}\ \ {\em $(i)$ If $\gamma>\gamma_*$,
  then the stationary solution $0$ of problem $(\ref{3.35})$ 
  is asymptotically stable. More precisely,  there exists 
  a positive constant $\epsilon$ such that for any given $\rho_0
  \in \mathcal O_\delta$ with
  $\|\rho_0\|_{h^{4+\alpha}(\Bbb S)}<\epsilon$, 
  problem $(\ref{3.35})$ has a unique solution 
  $\rho(t)\in C([0,+\infty),\mathcal O_\delta)
  \cap C^1([0, +\infty)$,
  $h^{1+\alpha}(\Bbb S))$, which converges exponentially fast
  to $0$ as $t\to +\infty$.

  $(ii)$ If $0<\gamma<\gamma_*$,
  then the stationary solution $0$ is unstable.
}

\medskip

{\bf Proof}. \ \ $(i)$ Let $\gamma >\gamma_*$.
  Recall that $h^{4+\alpha}(\Bbb S)$ is densely and compactly
  embedded into $h^{1+\alpha}(\Bbb S)$. Set $A:= -D \Psi(0)$
  and 
$$
  G(\rho):=-\Psi(\rho)+D\Psi(0)\rho \qquad\mbox{ for}\;\;
   \rho\in\mathcal O_\delta.
$$
  Clearly,
  we have $G(0)=0$ and $D G(0)=0$.
  Problem (\ref{3.35}) is equivalent to the following problem
\begin{equation}\label{5.1}
  \rho'(t)=A\rho(t)+G(\rho(t)) \qquad\mbox{for}\;\;t>0,\qquad\quad
  \rho(0)=\rho_0.
\end{equation}
  By Lemma 3.6, $A$ generates a strongly continuous 
  analytic semigroup on $h^{1+\alpha}(\Bbb S)$.  By Corollary 4.3 $(i)$,
  we have $\sup \{ \mbox {Re} \lambda: \lambda\in\sigma(A)\}
  <-\varpi<0$.  Thus by Theorem 9.1.2 of [\ref{lunardi}], 
  there are positive constants  $\omega, \epsilon$ and $M$ 
  such that if the initial value $\rho_0\in \mathcal{O}_\delta$ 
  and $\|\rho_0\|_{h^{4+\alpha}(\Bbb S)}<\epsilon$, then 
  the solution $\rho(t)$ 
  of problem (\ref{3.35}) exists globally and 
\begin{equation}\label{5.2}
  \|\rho(t)\|_{h^{4+\alpha}(\mathbb{S})}+\|\rho'(t)\|_{h^{1+\alpha}(\Bbb S)}
  \leq M e^{-\omega t}\|\rho_0\|_{h^{4+\alpha}
  (\mathbb{S})} \qquad \mbox {for}\;\; t\geq 0.
\end{equation}

  $(ii)$ If $0<\gamma<\gamma_*$, by Corollary 4.3 $(ii)$ we have 
 $\sigma_+(A)= \sigma(A)\cap \{\lambda\in \Bbb C: {\rm Re}\lambda>0\}\neq\emptyset$
 and $\inf \{ {\rm Re} \lambda: \lambda\in\sigma_+(A)\}>0$.
 Thus by Theorem 9.1.3 in [\ref{lunardi}], the stationary solution $\rho=0$ is unstable.
 The proof is complete.\qquad $\Box$

\medskip

  {\bf The proof of Theorem 1.2} \ \   By Lemma 3.1, Lemma 3.5 and Theorem 5.1 $(i)$,
  we see that  the flat stationary solution $(\sigma_s, p_s, \eta_s, \rho_s)$
  is asymptotically stable for $\gamma>\gamma_*$. More precisely,
  there is a constant $\epsilon>0$ such that for any
  $\rho_0\in \mathcal O_\delta$  satisfying
  $\|\rho_0\|_{h^{4+\alpha}(\Bbb S)}<\epsilon$,
  problem $(\ref{1.1})$ has a unique
  global solution $(\sigma(t), p(t), \eta(t), \rho(t))$ 
  with the form of
$$
  \sigma(t)=\Phi_*^{\tilde\rho(t)} u(t),\qquad
  p(t)=\Phi_*^{\tilde\rho(t)} v(t), \qquad \eta(t)=\eta_s+\mathcal S(\tilde\rho(t)),
  \qquad \rho(t)=\rho_s+\tilde\rho(t),
$$
  where $\tilde\rho(t)$ is the solution of problem (\ref{3.35}) with
  $\tilde\rho(0)=\rho_0$, and $u(t)$, $v(t)$, $\mathcal S(\tilde\rho(t))$
  are given by (\ref{3.17}) and (\ref{3.32}). By (\ref{5.2}) and the reduction in
  Section 3, we see that $(\sigma(t),p(t), \eta(t), \rho(t))$ 
  converges exponentially fast to  $(\sigma_s, p_s, \eta_s, \rho_s)$ 
  in $h^{4+\alpha}(\Omega_{\tilde\rho(t)}\backslash J_{\mathcal S(\tilde\rho(t))} )
  \times h^{2+\alpha}
  (\Omega_{\tilde\rho(t)}\backslash J_{\mathcal S(\tilde\rho(t))} )
  \times h^{4+\alpha}(\Bbb S) \times h^{4+\alpha}(\Bbb S)$,
  as time goes to infinity.

  Similarly,  by Lemma 3.1, Lemma 3.5 and
  Theorem 5.1 $(ii)$, 
  the flat stationary solution $(\sigma_s, p_s, \eta_s, \rho_s)$
  is unstable for $0<\gamma<\gamma_*$.  
  The proof is complete.
  \qquad $\Box$

\medskip

 {\bf Remark 5.2} \ \ From (\ref{4.19}) and (\ref{4.21}), we easily obtain 
  $d\gamma_k/d\nu<0$ for each $k\neq0$, $k\in \Bbb Z$. Thus we have
  $d\gamma_*/d\nu\le0$. It implies that the 
  smaller value of $\nu$ may make
  tumor more unaggressive.  In the limiting case $\nu=0$,
  since $\lambda_0(\gamma)= \nu=0$,  we have
  $0\in\sigma(D\Psi(0))$. It implies that the flat stationary 
  solution is not asymptotically stable any more for all $\gamma>0$.

\medskip
\hskip 1em


\vsss

\end{document}